\documentclass[12pt]{article}
\usepackage{jhep-mod-2019}
\usepackage[pdftex]{hyperref}
\hypersetup{
     colorlinks = true,
     linkcolor = blue,
     anchorcolor = blue,
     citecolor = blue,
     filecolor = blue,
     urlcolor = blue
     }
\usepackage{tikz}
\definecolor{purple}{rgb}{1,0,1}
\definecolor{lime}{HTML}{A6CE39} % needs xcolor

\newcommand{\blue}[1]{{\color{blue} #1}}

\newtheorem{theorem}{Theorem}

\definecolor{lime}{HTML}{A6CE39}
\newcommand{\orcidicon}{%
	\begin{tikzpicture}
	\draw[lime, fill=lime] (0,0) 
		circle [radius=0.16] 
		node[white] {{\fontfamily{qag}\selectfont \tiny ID}};
	\draw[white, fill=white] (-0.0625,0.095) 
		circle [radius=0.007];
	\end{tikzpicture}
	\hspace{-2mm}
}
\newcommand\orcidMatt{{\href{https://orcid.org/0000-0003-1088-6485}{\orcidicon}}}
\begin{document}
%========================================================
\title{\huge Strong version of Andrica's conjecture}

\author{\Large Matt Visser\,\orcidMatt{}}

\affiliation{
School of Mathematics and Statistics,
Victoria University of Wellington, \\ PO Box 600, Wellington 6140, 
New Zealand.
}

\emailAdd{matt.visser@sms.vuw.ac.nz}

\abstract{
\parindent0pt
\parskip7pt

A strong version of Andrica's conjecture can be formulated as follows:
Except for $p_n\in\{3,7,13,23,31,113\}$, that is $n\in\{2,4,6,9,11,30\}$, one has
$\sqrt{p_{n+1}}-\sqrt{p_n} < \frac{1}{2}.$
While a proof is far out of reach I shall show that this strong version of Andrica's conjecture is unconditionally and explicitly verified for all primes below the location of the 81$^{st}$ maximal prime gap, certainly for all primes $p <2^{64}\approx 1.844\times 10^{19}$. 
Furthermore this strong Andrica conjecture is slightly stronger than
Oppermann's conjecture --- which in turn is slightly stronger than both the strong and standard Legendre conjectures, and the strong and standard Brocard conjectures. 
Thus the Oppermann conjecture, and strong and standard Legendre conjectures, are  all unconditionally and explicitly verified for all primes $p <2^{64}\approx1.844\times 10^{19}$. 
Similarly, the strong and standard Brocard conjectures are  unconditionally and explicitly verified for all primes $p <2^{32} \approx 4.294 \times 10^9$. 

{\sc Sequences:} A005250 A002386 A005669 A000101 A107578 A050216 A050216 



{\sc arXiv:} 1812.02762   [math.NT]

{\sc Date:} 5 December 2018; 19 February 2019; 10 May 2019; \LaTeX-ed \today
}

%%%%%%%%%%%%%%%%%%%%%%%%%%%%%%%%%%%%%%%%%%
%%%%%%%%%%%%%%%%%%%%%%%%%%%%%%%%%%%%%%%%%%
%%%%%%%%%%%%%%%%%%%%%%%%%%%%%%%%%%%%%%%%%%
\maketitle
%========================================================
\section{Introduction}
%========================================================
%---------------------------------------------------------------------------------------------------------------------------------------------
\label{S:intro}
%---------------------------------------------------------------------------------------------------------------------------------------------
\def\N{{\mathbb{N}}}
\def\Z{{\mathbb{Z}}}
\def\P{{\mathbb{P}}}
\def\implies{\Longrightarrow}
\newtheorem{conjecture}{Conjecture}{}

The distribution of the primes is a fascinating topic that continues to provide many subtle and significant open questions~\cite{Ribenboim:91, Ribenboim:96, Ribenboim:04, Wells, Cramer:1919, Cramer:1936, Goldston, Rosser:38, Rosser:41, Cesaro, Cippola, Rosser:62, Sandor, Dusart-1999, Lowry-Duda, Dusart:2010, Trudigan:2014, Dusart2, Axler:2017, Visser:Lambert, Andrica, Visser:Andrica}.
In the current article I will consider strong, standard, and weak versions of the Andrica conjecture, and the closely related Oppermann, Legendre, and Brocard conjectures. These conjectures all impose constraints on the prime gaps of the form 
\begin{equation}
g_n := p_{n+1}-p_n = O\left( \sqrt{p_n}\,\right).
\end{equation}
Specifically, consider the strengthened version of the usual Andrica conjecture presented below.
\begin{conjecture}{(Strong Andrica conjecture)}\\
Except for $p_n\in\{3,7,13,23,31,113\}$, that is $n\in\{2,4,6,9,11,30\}$, one has
\begin{equation}
 \sqrt{p_{n+1}}-\sqrt{p_n} < {1\over2}; \qquad   \hbox{equivalently} \qquad    g_n := p_{n+1}-p_n < p_n^{1/2} +{1\over4}.
\end{equation}
\end{conjecture}
We shall soon see that this strong Andrica conjecture, like the usual Andrica conjecture~\cite{Andrica}, 
has the virtue that it can easily be verified on suitable intervals by inspecting a table of known maximal prime gaps~\cite{Visser:Andrica,g80,g75,gaps,Nicely}. 

The specific choice of the constant $1\over2$ in the left sub-equation of (1.2), which leads to the coefficient unity in front of the $ p_n^{1/2}$ in the right sub-equation of (1.2), was (with hindsight) carefully arranged to be just strong enough to imply the Oppermann conjecture below. The specific offset $1\over4$ in the right sub-equation is merely $\left(1\over2\right)^2$ and is not particularly important, more on this later.
In counterpoint, Oppermann's conjecture~\cite{Ribenboim:91, Ribenboim:96, Ribenboim:04, Wells, Oppermann} can be cast in either of the two equivalent forms given below.
\begin{conjecture}{(Oppermann conjecture~\cite{Ribenboim:91, Ribenboim:96, Ribenboim:04, Wells, Oppermann})}\\
(1) For all integers $m\geq 2$ there is at least one prime in each of the intervals
\begin{equation}
 \left(m(m-1), m^2\right);   \qquad\hbox{and} \qquad  \left(m^2, m(m+1)\right). 
\end{equation}
(2) For all integers $m\geq 1$ there is at least one prime in each of the intervals
\begin{equation}
 \left(m^2, m(m+1)\right) \qquad\hbox{and} \qquad  \left(m(m+1), (m+1)^2\right).
\end{equation}
\end{conjecture}
No proof of Oppermann's conjecture is known as of February 2019, so one must instead  resort to verifying it on certain (hopefully large) intervals. For this purpose it is useful to note that the strong variant of Andrica's conjecture introduced above implies Oppermann's conjecture.

Other weaker conjectures closely related to Oppermann's conjecture are:
\begin{conjecture}{(Strong Legendre conjecture)}\\
There are at least two primes in the interval
\begin{equation}
\left( m^2, (m+1)^2\right). 
\end{equation}
\end{conjecture}
\begin{conjecture}{(Standard Legendre conjecture~\emph{\cite{Legendre, Legendre-online, Landau-online})}}\\
There is at least one prime in the interval
\begin{equation}
\left( m^2, (m+1)^2\right). 
\end{equation}
\end{conjecture}
\begin{conjecture}{(Strong Brocard conjecture)}\\
There are at least  $2g_n:=2(p_{n+1}-p_n)$ primes in the interval
\begin{equation}
\left( p_n^2, p_{n+1}^2\right). 
\end{equation}
\end{conjecture}
\begin{conjecture}{Standard Brocard conjecture~\emph{\cite{Ribenboim:91, Ribenboim:96, Ribenboim:04, Brocard})}}\\
There are at least four primes in the interval
\begin{equation}
\left( p_n^2, p_{n+1}^2\right). 
\end{equation}
\end{conjecture}

For completeness we also define:
\begin{conjecture}{(Standard Andrica conjecture~\emph{\cite{Ribenboim:91, Ribenboim:96, Ribenboim:04, Wells, Andrica})}}\\
 Either of these two equivalent forms
\begin{equation}
\forall n \geq 1:\qquad
\sqrt{p_{n+1}}-\sqrt{p_n} <1; \qquad g_n:= p_{n+1}-p_n < 2\sqrt{p_n}+1.
\end{equation}
\end{conjecture}
\begin{conjecture}{(Various weakened versions of the Andrica conjecture)}\\
A weakened version of the Andrica conjecture can be presented in either of these two equivalent forms
\begin{equation}
\forall n \geq 1:\qquad
\sqrt{p_{n+1}}-\sqrt{p_n} <2; \qquad equivalently \qquad g_n:= p_{n+1}-p_n < 4\sqrt{p_n}+4.
\end{equation}
We could try to be even more general and conjecture
\begin{equation}
\forall n \geq 1, \; \forall c > 1:\quad equivalently \quad
\sqrt{p_{n+1}}-\sqrt{p_n} <c; \quad g_n:= p_{n+1}-p_n < 2c\sqrt{p_n}+c^2,
\end{equation}
but the specific choice $c=2$ is more useful in that we shall soon see that it is easily related to the standard Legendre conjecture.
An even weaker conjecture would be to merely assert that $\sqrt{p_{n+1}}-\sqrt{p_n}$ is bounded. 
\end{conjecture}

\vspace{-10pt}
%==================================================================
\section[Verifying the strong Andrica conjecture for primes $p < 2^{64}$]
{{Verifying the strong Andrica conjecture for primes $p < 2^{64}$}}
%==================================================================
%---------------------------------------------------------------------------------------------------------------------------------------------
\label{S:verify-gap}
%---------------------------------------------------------------------------------------------------------------------------------------------
The argument is a minor variant of that given for the standard Andrica conjecture in reference~\cite{Visser:Andrica}.
Consider the maximal prime gaps: Following the notation developed in reference~\cite{Visser:Andrica}, let the triplet $(i, g^*_i, p^*_i)$ denote the $i^{th}$ maximal prime gap; of width $g^*_i$, starting at the prime $p^*_i$. 
(See see the sequences A005250, A002386, A005669, A000101, A107578.)
80 such maximal prime gaps are currently known~\cite{g80,g75,gaps}, up to $g^*_{80}=1550$ and $p^*_{80}= 18,361,375,334,787,046,697 >1.836\times 10^{19}$.
Now consider the interval $[p^*_i, p^*_{i+1}-1]$, from the lower end of the $i^{th}$ maximal prime gap to just below the beginning of the $(i+1)^{th}$ maximal prime gap. Then everywhere in this interval
\begin{equation}
\forall p_n\in [p^*_i,p^*_{i+1}-1]  \qquad   g_n \leq g^*_i; \qquad  \sqrt{p_i^*} +{1\over4} \leq  \sqrt{p_n} +{1\over4}.
\end{equation}
\enlargethispage{25pt}
Therefore, if the strong Andrica conjecture holds at the beginning of the  interval $p_n\in [p^*_i,p^*_{i+1}-1]$, then it  certainly holds on the entire interval. Explicitly checking a table of maximal prime gaps~\cite{g80}, the strong Andrica conjecture certainly holds on the interval $[p_7^*,p^*_{81}-1]$, from $p^*_7=523$  up to just before the beginning of the 81$^{st}$ maximal prime gap, $p^*_{81}-1$, even if we do not yet know the value of $p^*_{81}$. Then explicitly checking the primes below $p_7^*= 523$ the strong Andrica conjecture holds for all primes $p$ less than $p_{81}$ except $p\in\{3,7,13,23,31,113\}$. 

Even though we do not explicitly know $p^*_{81}$ a safe fully explicit statement is this:  Use the bound $p^*_{81} > 2^{64}$ that comes from an exhaustive search for all maximal prime gaps below $2^{64}$ (see reference~\cite{Nicely}) to observationally verify 
 the strong Andrica conjecture for all primes $p < 2^{64} = 18,446,744,073,709,551,616 \approx 1.844\times 10^{19}$. (This argument also automatically verifies the standard and weak Andrica conjectures over the same domain.)

%========================================================
\section
[Verifying the Oppermann conjecture for primes $p<2^{64}$\\ and integers $m< 2^{32}$]
{Verifying the Oppermann conjecture for primes $p<2^{64}$ \\ and integers $m<2^{32}$}
%========================================================
%---------------------------------------------------------------------------------------------------------------------------------------------
\label{S:verify-oppermann}
%---------------------------------------------------------------------------------------------------------------------------------------------

\begin{theorem}
The strong Andrica conjecture implies the Oppermann conjecture.
\end{theorem}
Proof:\\
(1) Note $\left[\sqrt{113}\,\right] = 10$, so taking $m\geq12$ means we safely avoid the exceptional cases  in the strong Andrica conjecture.  
\\
(2) Pick some fixed $m\geq 12$ and let $p_n$ be the largest prime less than $m^2$. \\
Then by construction $p_{n}< m^2<p_{n+1} $ and by the strong Andrica conjecture we have
\begin{equation}
p_{n+1} := p_n + g_n  < p_n  + p_n^{1/2} +{1\over4} <  m^2 + m + {1\over4}.
\end{equation}
But since $p_{n+1}\in \N$ this implies $p_{n+1} \leq m(m+1)$.\\
But since $p_{n+1}\in \P$ this implies $p_{n+1} < m(m+1)$.\\
(3) Pick some fixed $m\geq 12$ and let $p_{n+1}$ be the smallest prime greater than $m^2$. \\
Then by construction $p_{n+1}> m^2> p_n$ and by the strong Andrica conjecture we have
\begin{equation}
p_{n} := p_{n+1}-  g_n  > p_{n+1}  - p_n^{1/2}  - {1\over 4} > m^2  - m  - {1\over 4}.
\end{equation}
But since $p_{n}\in \N$ this implies $p_{n} \geq m(m+1)$.\\
But since $p_{n}\in \P$ this implies $p_{n} > m(m+1)$.\\
(4) Check the cases $m<12$ by explicit computation.

Now while the proof that the strong Andrica conjecture implies the Oppermann conjecture is unconditional, we have only verified the strong Andrica conjecture up to the location of the 81$^{st}$ maximal prime gap $p< p^*_{81}$. Consequently we can only verify the Oppermann conjecture up to $m \leq \left[\sqrt{p^*_{81}}\,\right]$. Since we do not yet know $p^*_{81}$ the best explicit verification range is based on the bound $p^*_{81} > 2^{64}$ (see reference~\cite{Nicely}) implying explicit verification for integers
$m \leq 2^{32} = 4,294,967,296\approx 4.294 \times 10^9$.

%========================================================
\section{Verifying other weaker conjectures}
%========================================================
%---------------------------------------------------------------------------------------------------------------------------------------------
\label{S:verify-weaker}
%---------------------------------------------------------------------------------------------------------------------------------------------
\vspace{-15pt}
First, note that the Oppermann conjecture implies the strong Legendre conjecture which in turn implies the standard Legendre conjecture.
To see this write
\begin{equation}
(m^2, (m+1)^2) =  (m^2, m(m+1)) \cup \{m(m+1)\} \cup (m(m+1), (m+1)^2).
\end{equation}
Then assuming the  Oppermann conjecture, the interval $ (m^2, m(m+1))$ contains a prime. Also, the interval $ (m(m+1), (m+1)^2)$ contains a prime,
and $m(m+1)$ is not a prime. So there are at least two primes in $(m^2,(m+1)^2)$. Furthermore,  automatically there is at least one prime in $(m^2,(m+1)^2)$.

Second, note that the strong Legendre conjecture implies the strong Brocard conjecture, which in turn implies the standard Brocard conjecture. To see this, split the interval $(p_n^2,p_{n+1}^2)$ into $g_n:= p_{n+1}-p_n$ sub-intervals of the form $\left( (p_n+i)^2,(p_n+i+1)^2\right)$ for $0\leq i \leq g_n-1$. Assuming the strong Legendre conjecture, each of these $g_n$ sub-intervals contains at least 2 primes, so $(p_n^2,p_{n+1}^2)$ contains at least $2g_n$ primes, which is the strong Brocard conjecture. Since $g_n\geq 2$ this automatically implies the standard Brocard conjecture.

Third, the Oppermann conjecture implies the standard Andrica conjecture.
To see this define the two integers $m=[p_{n}^{1/2}]$ and  $\tilde m=[p_{n+1}^{1/2}]$. 
\begin{itemize}
\item 
If $m=\tilde m$ then $p_n > m^2$ and $p_{n+1} < (m+1)^2$.
So $\sqrt{p_{n+1}}-\sqrt{p_n} < (m+1)-m = 1$.
\item 
If $m\neq \tilde m$ then $ m^2 < p_n < \tilde m^2 < p_{n+1}$. Then, assuming the Oppermann conjecture, we have both $p_n >  \tilde m^2-\tilde m$ and $p_{n+1} < \tilde m^2 + \tilde m$, and so $g_n < 2\tilde m $.
But then 
\begin{equation}
p_n >  \tilde m(\tilde m-1) = (\tilde m-1/2)^2-1/4 > (\tilde m-1/2)^2,
\end{equation}
so $ \tilde m < 1/2 + \sqrt{p_n}$, and so $2\tilde m< 2\sqrt{p_n}+1$.
Finally this implies $g_n < 2\sqrt{p_n}+1$. 
\end{itemize}

\enlargethispage{30pt}
Fourth, the standard  Andrica conjecture implies the standard Legendre conjecture. 
Pick some fixed $m$ and let $p_n$ be the largest prime less than $m^2$. Then by construction $p_{n+1}> m^2> p_n$ and assuming the standard Andrica conjecture we have
\begin{equation}
p_{n+1} = p_n + g_n  < p_n  + 2 \sqrt{p_n} +1 <  m^2 + 2m + 1 = (m+1)^2.
\end{equation}
Fifth, the standard  Legendre conjecture implies the weak $(c=2)$ Andrica conjecture. To see this pick some prime $p_n$ and let $m=[\sqrt{p_n}]$, then $m^2 < p_n < (m+1)^2$ and by the standard Legendre conjecture there is at least one more prime in $\left((m+1)^2, (m+2)^2\right)$. Then 
\begin{equation}
g_n := p_{n+1} - p_n < (m+2)^2-m^2 = 4 m + 4 < 4 \sqrt{p_n} + 4,
\end{equation}
which is the weak Andrica conjecture.

In view of what we have already seen for the Oppermann conjecture, we now see that the strong and weak Legendre conjectures
are certainly verified for integers up to $m \leq \left[\sqrt{p^*_{81}}\,\right]$. Since we do not yet know $p^*_{81}$ the best explicit verification range comes from the bound $p^*_{81}< 2^{64}$ (see reference~\cite{Nicely}) corresponding to integers $m \leq 2^{32} \approx 4.294 \times 10^9$.
Similarly, the strong and weak Brocard conjectures are certainly verified for primes up to $p \leq \left[\sqrt{p^*_{81}}\,\right]$. Since we do not yet know $p^*_{81}$ the best explicit verification range is for primes $p < 2^{32} \approx 4.294 \times 10^9$.

%------------------------------------------------
\section{Discussion}\label{S:Discussion}
%------------------------------------------------
While it is reasonably well known that all of the the Andrica, Oppermann,  Legendre, and Brocard conjectures are closely related, little work seems to have gone into using these inter-relations to find suitably large regions where these conjectures can all be verified to be true. By setting up a strong version of the Andrica conjecture one can easily demonstrate that all of the  strong, standard, and weak Andrica conjectures are certainly valid for primes $p < 2^{64}$.
By proving that the strong Andrica conjecture implies the Oppermann conjecture, which in turn implies the strong and standard Legendre conjectures, and the strong and standard Brocard conjectures, one can demonstrate that the Oppermann, and strong and standard Legendre conjectures, are likewise certainly all valid for primes $p < 2^{64}$, corresponding to integers $m \leq 2^{32}$. Similarly the strong and weak Brocard conjectures are certainly valid for 
primes  $ p < 2^{32}$.
The key item in these bounds is the location of the highest-known maximal prime gap, so updating these bounds will be automatic as new maximal prime gaps are identified. 

\enlargethispage{35pt}
In counterpoint, what would it take for all of these conjectures to suddenly fail at the next opportunity, the 81$^{st}$ maximal prime gap? One would need $ g_{81}^* > \sqrt{p_{81}^*}+{1\over 4} >\sqrt{p_{80}^*} \sim 4.285 \times 10^9$. That is, since $g_{80}^* = 1550$, one would need the next maximal prime gap to suddenly change from order $10^3$ to order $10^9$. While, given current knowledge, this cannot be entirely ruled out --- it does at the very least look suggestively unlikely.

Finally, while the strong, standard, and weak versions of the Andrica conjecture, (and the closely related Oppermann, Legendre, and Brocard conjectures) all impose constraints on the prime gaps of the form 
\begin{equation}
g_n := p_{n+1}-p_n = O\left( \sqrt{p_n}\,\right),
\end{equation}
these are by no means the most stringent conjectures one might plausibly make. 
For instance, consider the following. 
\begin{conjecture}{(Square root conjecture)}\\
Except for $p_n\in\{3,7,13,23,31,113\}$, that is $n\in\{2,4,6,9,11,30\}$, one has
\begin{equation}
g_n := p_{n+1}-p_n < p_n^{1/2}. 
\end{equation}
\end{conjecture}
This is equivalent to asserting
\begin{equation}
\sqrt{p_{n+1}}-\sqrt{p_n} < \sqrt{p_n} \left(\sqrt{1+{1\over\sqrt{p_n}}}-1\right). 
\end{equation}

This square root conjecture can also be easily verified to certainly hold for all primes less than $p_{81}^*$. The price paid here is that while the conjecture looks somewhat simpler when phrased in terms of prime gaps, (no $1\over4$ offset term), the statement in terms of $\sqrt{p_{n+1}}-\sqrt{p_n}$ is more complicated and less ``Andrica-like'' in flavour. Note that 
\begin{equation}
\sqrt{p_n} \left(\sqrt{1+{1\over\sqrt{p_n}}}-1\right) < {1\over2}, 
\end{equation}
indeed
\begin{equation}
 \sqrt{p_n} \left(\sqrt{1+{1\over\sqrt{p_n}}}-1\right)= {1\over2} - {1\over 8\sqrt{p_n}} + O\left(1\over p_n\right), 
\end{equation}
so this square root conjecture asymptotically approaches the strong Andrica conjecture for large primes. 
The current best unconditional result along these lines is the Baker--Harman--Pintz~\cite{BHP} result that for sufficiently large $x$ the interval $[x-x^{0.525}, x]$ always contains primes --- so that for sufficiently large primes  $p_{n+1}-p_{n} \leq O\left((p_n)^{0.525}\right)$.
Note the exponent has not yet been (unconditionally) reduced to $1\over2$, and the implied constant in the phrase ``sufficiently large'' is still undetermined.

\clearpage
Furthermore, observe that in references~\cite{Ribenboim:91, Ribenboim:96, Ribenboim:04} the author mentions the ``open problem'' as to whether
\begin{equation}
\lim_{n \to \infty} \left( \sqrt{p_{n+1}} - \sqrt{p_{n}} \right)= 0? 
\end{equation}
If this limit exists, (and it is easy to see that $\liminf_{n \to \infty} \left( \sqrt{p_{n+1}} - \sqrt{p_{n}} \right)= 0$, this is for instance a special case of  the discussion in reference~\cite{Andrica}, so it is only the existence of the limit that is in question), then for any specified $\epsilon>0$ the inequality $ \sqrt{p_{n+1}} - \sqrt{p_{n}} < \epsilon$ can be violated only a finite number of times. 
This observation can be linked back to the finite ``exception list'' we needed to invoke in setting up the strong Andrica conjecture.

Finally the Cram\'er conjecture~\cite{Cramer:1936} (and closely related conjectures such as the Firoozbakht conjecture~\cite{Ribenboim:91, Ribenboim:96, Ribenboim:04,Firoozbakht,Firoozbakht2,Firoozbakht:2015a,Firoozbakht:2015b,Kourbatov3,Nicholson3})
impose significantly stronger constraints on the prime gaps of the form
\begin{equation}
g_n := p_{n+1}-p_n = O\left( (\ln p_n)^2\right).
\end{equation}
I will not say anything further regarding the Cram\'er and Firoozbakht conjectures in the current article, but hope to turn to this topic in future work.

%------------------------------------------------
\acknowledgments{
This research was supported by the Marsden Fund, administered by the Royal Society of New Zealand.  
I particularly wish to thank Alexei Kourbatov for his input regarding known bounds on the location of the 81$^{st}$ maximal prime gap ($p^*_{81}>2^{64}$). 
}
%-------------------------------------------------

%\conflictsofinterest{The author declares no conflict of interest.}

%\clearpage

\vspace{-10pt}

%\reftitle{References and Notes}
%========================================================

%========================================================
\end{document}